\documentclass[3p,,preprint,12pt]{elsarticle}

\usepackage{graphicx} 
\usepackage{amsmath}
\usepackage{amssymb}
\usepackage{amsthm}
\usepackage{lineno}
\usepackage{color}
\usepackage{verbatim}
\usepackage{subcaption}
\usepackage{url,multirow,morefloats,floatflt,cancel,tfrupee}
\makeatletter
\let\save@ps@pprintTitle\ps@pprintTitle
\def\ps@pprintTitle{\save@ps@pprintTitle\gdef\@oddfoot{\footnotesize\itshape \null\hfill\today}}
\makeatother

\newtheorem{theorem}{Theorem}

\newtheorem{proposition}{Proposition}
\newtheorem{definition}{Definition}
\newtheorem{corollary}{Corollary}

\def\F{{\cal F}}
\def\L{{\cal L}}
\def\M{{\cal M}}
\def\sumi{\sum_{i=1}^n}

\journal{Operations Research Letters}

\begin{document}
\begin{frontmatter}

\title{On monotone completion of risk markets: Limit results for incomplete risk markets}

\author[UMass]{Iman Khajepour\corref{cor1}}
\ead{ikhajepourta@umass.edu}  
\author[Auckland]{Geoffrey Pritchard}
\ead{g.pritchard@auckland.ac.nz}  
\author[JBS]{Danny Ralph}
\ead{d.ralph@jbs.cam.ac.uk} 
\author[UMass]{Golbon Zakeri}
\ead{gzakeri@umass.edu}  

\affiliation[UMass]{%
            organization={Department of Industrial Engineering}, 
            addressline={University of Massachusetts Amherst}, 
            city={Amherst}, postcode={MA 01003}, state={MA}, country={USA}}
\affiliation[Auckland]{%
            organization={Department of Statistics}, 
            addressline={University of Auckland}, 
            city={Auckland}, postcode={Private Bag 92019}, country={New Zealand}}
\affiliation[JBS]{%
            organization={Judge Business School}, 
            addressline={University of Cambridge}, 
            city={Cambridge}, postcode={CB2 1AG}, country={United Kingdom}}

\cortext[cor1]{Corresponding author}

\begin{abstract}
We consider a competitive market with risk-averse participants.
We assume that agents' risks are measured by coherent risk measures introduced by Artzner et al. (1999) \cite{Artzneretal}. Fundamental theorems of welfare economics have long established the equivalence of competitive equilibria and system welfare optimization (see e.g. Samuelson (1947) \cite{SamuelsonFEA}). 
These have been extended to the case of 
risk-averse agents with complete risk markets in Ralph and Smeers (2015) \cite{Ralph-Smeers-2015}. In this paper we consider risk trading in incomplete markets, and introduce a mechanism to complete the market iteratively while monotonically enhancing welfare. 
\end{abstract}
\end{frontmatter}

\section{Introduction and background}
\label{Introduction and background}
Over the past few decades, there has been a significant paradigm shift from centralized governance of various goods and services. Industries such as telecommunications, airlines, trucking, and energy have experienced substantial deregulation in the last half-century. Proponents of market systems argue that they better incentivize innovation, allocate resources more efficiently, provide greater choices, and offer increased flexibility compared to central planning.

Economists have extensively examined the comparative efficiency of these systems, with studies consistently showing lower efficiency in planned economies (see e.g., \cite{BergsonAER, BergsonJEP, MoroneyJCE}). However, an essential attribute of any resource allocation mechanism—whether market-based or centrally planned—is its ability to maximize overall societal benefit or welfare. 
Blaug \cite{BlaugHPE} gives an economic history of the idea that markets can, under perfect competition, maximize social welfare and how this idea relates to the first fundamental theorem of welfare economics in which economic efficiency is demonstrated in terms of Pareto optimal allocations of commodities to individuals rather than total welfare of the system. 
Arrow \cite{ArrowProc51} showed the equivalence between Pareto optimal allocation and welfare maximization via convex analysis,\footnote{\cite{ArrowProc51} assumes strictly concave utility functions and convex compact feasible sets.} where the equilibrium price is that price at which demand equals supply.
This applies to markets in risk rather than traditional commodities; risk markets are indeed the subject of Arrow~\cite{ArrowRES54} on the trading of securities for money. {\footnote{\cite{ArrowRES54} does not go to the trouble of considering constraints, but this turns out to be without loss of generality given~\cite{ArrowProc51}.}}
This view of risk markets was reiterated by Ralph and Smeers \cite{Ralph-Smeers-2015} who analysed risk markets in which financial products are traded to maximize the welfare of the system and the prices of financial products emerge as dual variables once the optimal set of trades is known. In this setting, it is immediate that the optimal trades will yield a system-optimal allocation of risk, whether or not the available securities cover all future uncertainties.
While \cite{Ralph-Smeers-2015} focusses on complete markets, which is the situation in which all future uncertainties, coded as Arrow-Debreu securities, are tradeable, section~4 of that paper also looks at the incomplete case as an optimal welfare problem.

 In practice, most risk markets are incomplete, i.e., some risks or future contingencies are uninsurable no matter what the combination of tradeable securities.
 A well-known consequence of incompleteness is that for any perfectly competitive market equilibrium, the welfare of any agent can be decreased (or increased) by introducing a single new asset; see Hart \cite{HartJET} for seminal work on this question and Elul \cite{ElulJET} for a more developed view. While this classical approach focusses on Pareto optimality, it is intuitive that under welfare maximization, which is what a risk market achieves, adding any additional asset can only improve the total welfare. 
When introducing risk-trading instruments to complete the market, it is important to ensure that welfare is enhanced as more instruments become available. To our knowledge, this paper presents the first method for iterative market completion that consistently increases total welfare.  

In the rest of this paper, we will briefly review the set up and results of \cite{Ralph-Smeers-2015}, layout our market completion approach in the face of a stochastic future with finite discrete uncertainty, provide examples and establish our monotonic welfare enhancement result in Section \ref{Incomplete Markets of Bundled Trade Instruments: finite discrete case}. In Section \ref{Incomplete Markets of Bundled Instruments: discrete distributions with countably infinite states }, we extend the models by incorporating an infinite number of scenarios. We will then present a theorem to explore the asymptotic behavior of social welfare when a market completion approach, similar to the one introduced in Section \ref{Incomplete Markets of Bundled Trade Instruments: finite discrete case}, is applied. In Section \ref{Incomplete Markets of Bundled Trade Instruments: arbitrary probability space}, we extend our model to address continuous distribution of future outcomes where the set up requires more careful treatment using probability theory. We will present an example and establish very similar monotonic welfare enhancement results, as well as a risk price convergence result. 

\section{Incomplete Markets of Bundled Trade Instruments: finite discrete case}
\label{Incomplete Markets of Bundled Trade Instruments: finite discrete case}
Consider a competitive market with $I$ agents. In this this section, we will confine our attention to the case where there are a finite number of future scenarios, indexed by $s \in S$. Each agent $i \in I$ is endowed with a loss vector denoted by $z_i \in \mathbb{R}^{|S|}$, where each element represents the agent's loss in the corresponding scenario. Agents evaluate their risk using a coherent risk measure $r_i(z_i)$ (in particular $r_i$ is convex) \cite{Artzneretal}. We introduce the set $J$ that partitions the set of scenarios. We refer to elements $S_j \in J$ as scenario bundles and note that $\cup_j S_j = S$ and that $S_j \cap S_i = \emptyset$ for $i\ne j$. We now define a generalization of an Arrow-Debreu security using $S_j$. 
\begin{definition} 
Any security bundle $S_j \in J$, will pay out one unit in the realization of {\em any} scenario $s \in S_j$ and no payout in other scenarios. 
We refer to this set up as the incomplete market ${\mathcal{M}_J}$. 
\end{definition}

In our market set up ${\mathcal{M}_J}$, agent $i$ trades based on their trade vector $t_i \in \mathbb{R}^{|J|}$ to improve their position after participating in the risk market. Here, $t_{ij}$ represents the volume of bundle $S_j$ bought (positive value) or sold (negative value) or by agent $i$. Note that purchasing a volume $t_{ij}$ of bundle $S_j$ is equivalent to purchasing the 
volume $t_{ij}$ of each of the Arrow-Debreu securities contained in the bundle. Hence the payoff is a vector $w \in \mathbb{R}^{|S|}$ with an associated one dollar payment for any component $s$ in the ``realized" bundle $S_j$. We formulate the risk adjusted loss minimization problem of each agent below, suppressing the agent's index.
\begin{align}
    \text{min}_{w , \ t} \quad &\lambda^T t + r(z - w)  & \quad  & \label{eq:Agent Optimization} \\
    \text{s.t.} \ \ & w_{k} = t_{j}   &\forall k \in S_j   \quad & [\gamma_{k}] \label{eq:bundle_constraint}
\end{align}
Note that there are $|J|$ groups of constraints of type \ref{eq:bundle_constraint}, one for each $S_j \in J$. 
Additionally, the price vector $\lambda$ is determined such that in equilibrium, the trades are balanced, i.e.
\[\sum_i t_i = 0. \]

Similar to \cite{Ralph-Smeers-2015}, we can consider a single optimization problem corresponding to total social welfare maximization, 
subject to the trade of the available instruments (i.e. available bundles of Arrow-Debreu securities $S_j$). 
\begin{align}
    \min_{w,\ t} \quad & \sum_{i=1}^N r_i(z_i - w_i)  & \quad  & \label{eq:Social Optimization} \\
    \text{s.t.} \ \ & w_{ik} = t_{ij}   &{\text{for each $i$, }} \forall k \in S_j   \quad & [\gamma_{ik}] \nonumber \\
    & \sum_{i=1}^N t_i = 0 &  \quad & [\lambda] \nonumber
\end{align}

The market equilibrium for the incomplete market, where agents trade bundles of securities instead of Arrow-Debreu securities, is equivalent to the social optimization problem in (\ref{eq:Social Optimization}). This can be easily observed through the equivalence of the necessary and sufficient optimality conditions of the
(\ref{eq:Social Optimization}) to those of the agents' problems, appended with market clearing in the risk trade market. 

\subsection{Refining Bundles}
Now consider the case when the partition $J$ is refined. This means that each subset $S_j$ is split into two (or more; we assume two sub-bundles for exposition purposes and without loss of generality) further sets $S_{j_1}$ and $S_{j_2}$.  As a result of this refinement, the group constraints associated with $j$ for each agent’s optimization problem are relaxed. That is 
$$
\begin{array}{ll}
w_{ik} = t_{ij} & \forall k \in S_j
\end{array}
$$
will now be replaced with the more relaxed system 
$$
\begin{array}{ll}
w_{ik} = t_{ij_1} & \forall k \in S_{j_1} \\
w_{ik} = t_{ij_2} & \forall k \in S_{j_2}
\end{array}
$$
This relaxation, ensures that the social optimization problem's objective will improve as we refine the security bundles. 
Together with the equilibrium equivalence lemma from the above, we obtain the theorem below. 
\begin{theorem}
    Consider the incomplete market set up ${\mathcal{M}_J}$ where $J$ is a partition of the scenario set $S$. Let $J_1 \cup J_2 =J $ be a refinement of the partition $J$, that is for each $S_j \in J$, $S_j = S_{j_1} \cup S_{j_2}$, where $S_{j_1} \cap S_{j_2} = \emptyset$. Then the incomplete market ${\mathcal{M}_{{J_1} \cup {J_2}}}$ is an enhancement of 
    ${\mathcal{M}_J}$, in the sense that there exists an equilibrium, where the total welfare of the system is improved.  
\end{theorem}
\begin{proof}[proof]
    This simply follows from the equivalence of the equilibrium and total welfare optimization problems in \ref{eq:Agent Optimization} and \ref{eq:Social Optimization}, and the fact that refinement of $J$ into $J_1 \cup J_2$, as defined in the hypothesis is a relaxation of the \ref{eq:bundle_constraint} in the total social welfare optimization. \qedhere
\end{proof}

\medskip
{\bf Example 1.}
Consider a competitive market with two agents and 256 uniformly distributed possible scenarios. The agents are endowed with loss vectors $z_1$ and $z_2$ defined as follows.  
$$z_{1s} = \frac{2s}{256}, \forall s \in \{1,\ldots, 256\}$$ while
$$
z_{2s} = \left\{
\begin{array}{ll}
\frac{s}{256} & \text{for } s \in \{1,\ldots, 128\} \\
1-\frac{s-129}{256} & \text{otherwise}
\end{array}
\right . 
$$

\begin{figure}[htp!]
  \centering
  \includegraphics[width=0.65\textwidth]{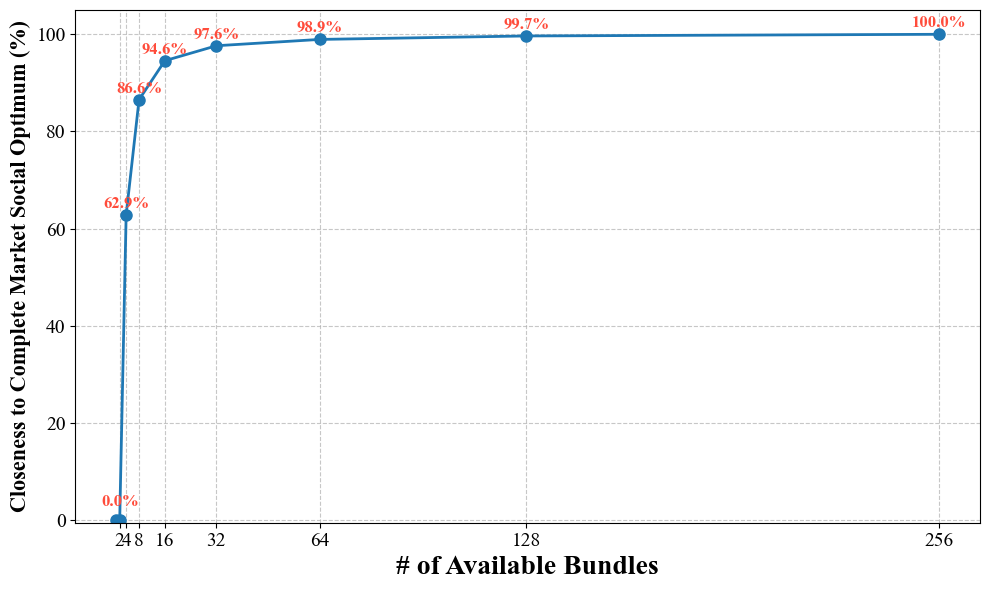} 
  \caption{\textbf{Improvement of social welfare through bundle refinement.} This figure illustrates the progression of incomplete markets towards the complete market social optimum and how refinement results in a non-decreasing trend.}
  \label{fig:exampleWelfareEnhancement}
\end{figure}

In this example, agents assess their risk using $CVaR_\beta$ (see e.g. \cite{Rock-Uryasev-2000}) with $\beta_1 = \frac{1}{5}$ and $\beta_2 = \frac{1}{4}$.
In each iteration $m \in \{0,\ldots,8\}$, the scenario set $S = \{1,\ldots,256\}$ is partitioned into $2^m$ bundles (subsets). Specifically, the first $2^{8-m}$ scenarios are grouped into the first bundle, the second $2^{8-m}$ scenarios into the second bundle, and this process continues until $2^m$ bundles are formed, each containing $2^{8-m}$ scenarios. During each iteration $m$, agents have access to these bundles for trading and risk mitigation.
To determine the equilibrium for the incomplete market and establish the prices of bundles, we solve the social optimization problem using the bundling approach for each specified value of $m$. Figure \ref{fig:exampleWelfareEnhancement} illustrates the enhancement of the system welfare through the refinement of bundles into smaller ones, assessing the difference between the incomplete market's optimal objective (i.e. total welfare) to that of the complete market's. 
We can also examine how bundle prices evolve as we refine the available bundles. The plots below visualize the prices of bundles, illustrating cases with 2, 4, and 8 bundles. 
\begin{figure}[htp!]
    \centering
    \begin{subfigure}{0.45\textwidth}
        \centering
        \includegraphics[width=\linewidth]{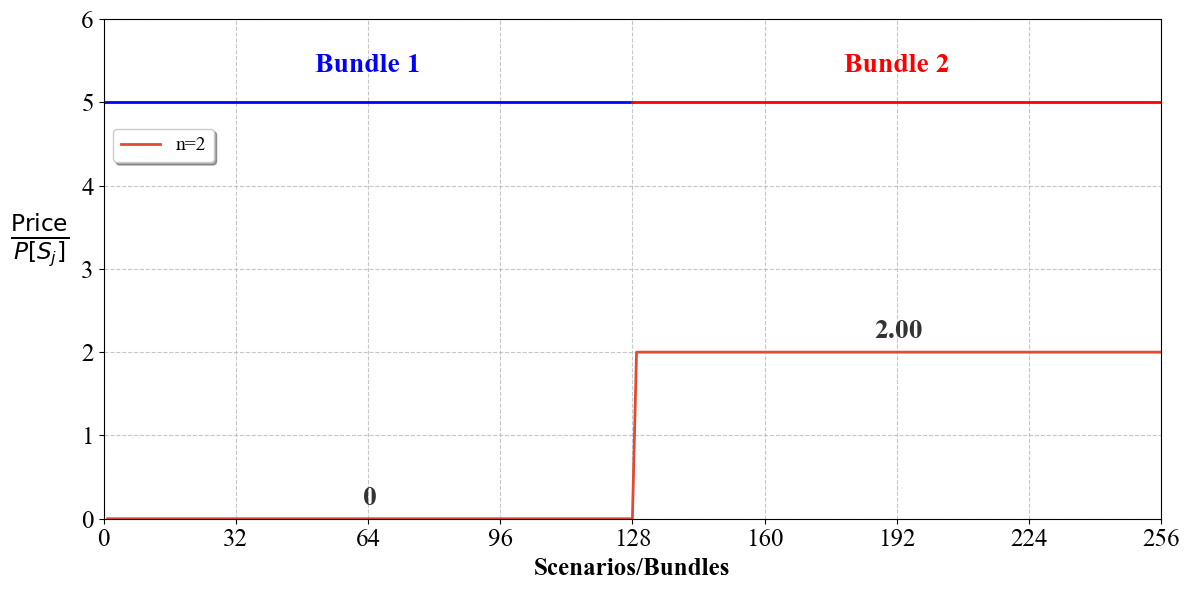} 
        \caption{Bundle prices for two available trade bundles.}
        \label{fig:img1}
    \end{subfigure}
    \hspace{0.05\textwidth} 
    \begin{subfigure}{0.45\textwidth}
        \centering
        \includegraphics[width=\linewidth]{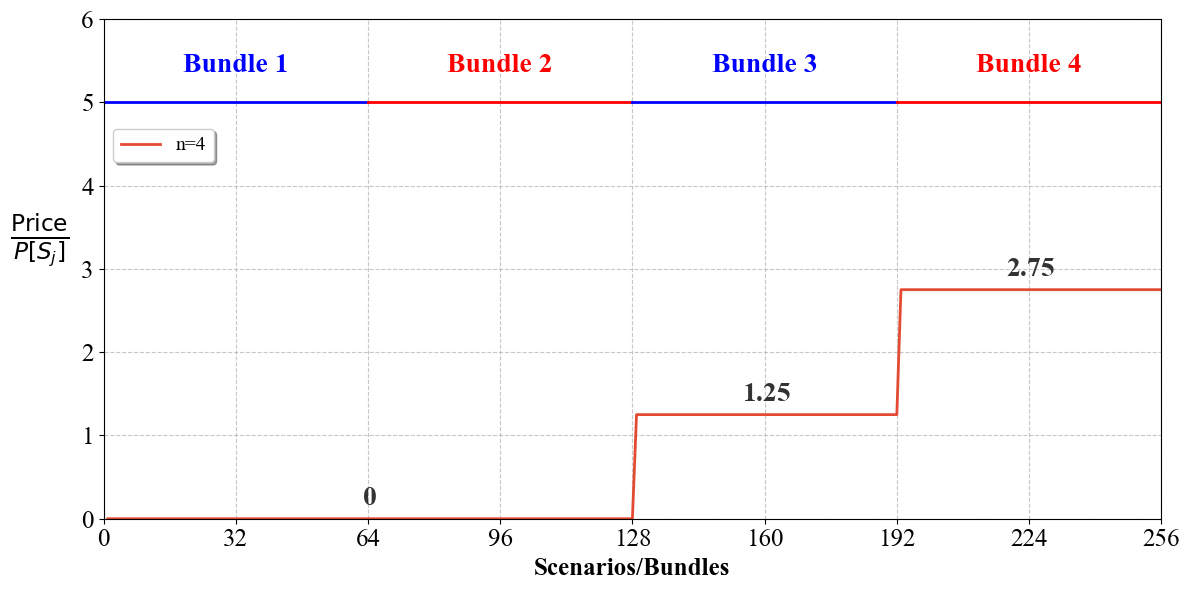} 
        \caption{Bundle prices for four available trade bundles.}
        \label{fig:img2}
    \end{subfigure}
    \vspace{0.05cm} 
    \begin{subfigure}{0.45\textwidth}
        \centering
        \includegraphics[width=\linewidth]{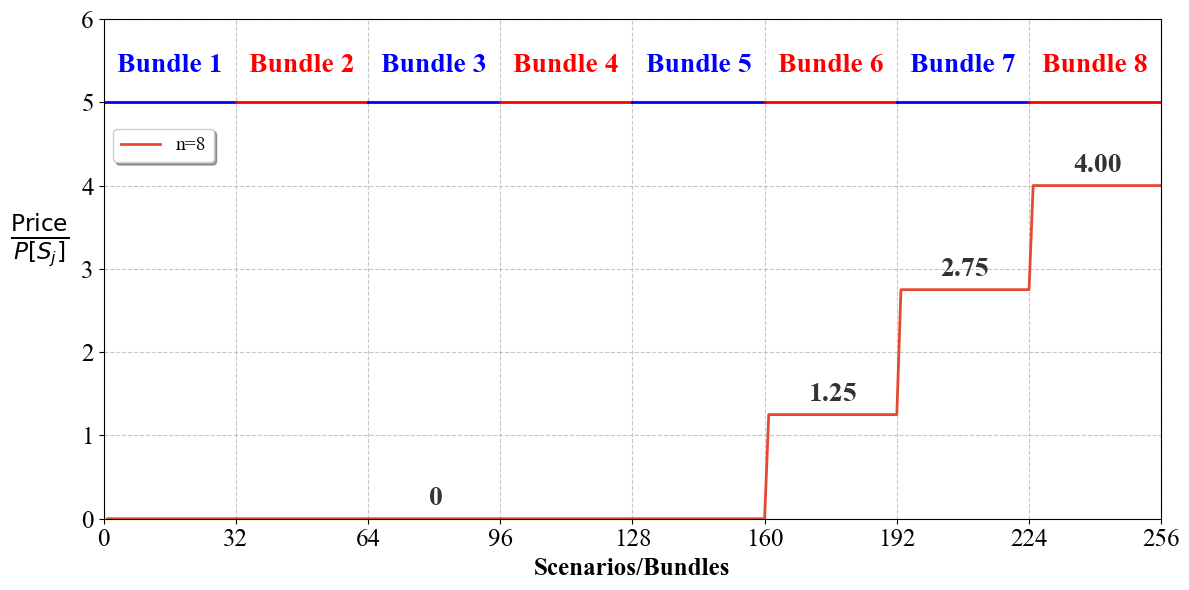} 
        \caption{Bundle prices for eight available trade bundles.}
        \label{fig:img3}
    \end{subfigure}
    \caption{Evolution of risk prices as bundles refine.}
    \label{fig:all_images}
\end{figure}
Continuing this process of bundle refinement for $m= 4 , \ldots ,8$ it becomes evident that as the market approaches completeness, the prices of bundles converge to a unique price, which corresponds to the price set by the least risk-averse agent in the complete market. Figure \ref{fig:five_images} illustrates this phenomenon.
\begin{figure}[htp!]
    \centering
    \begin{subfigure}{0.30\textwidth}
        \centering
        \includegraphics[width=\linewidth]{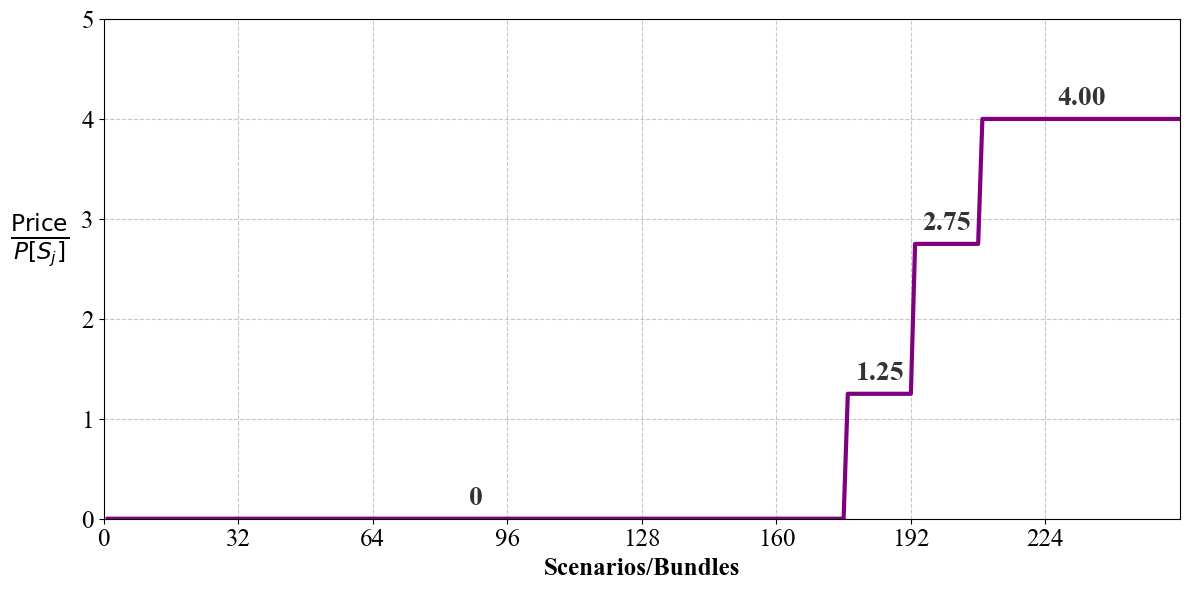}
        \caption{Bundle prices at iteration $m = 4$.}
    \end{subfigure}
    \begin{subfigure}{0.30\textwidth}
        \centering
        \includegraphics[width=\linewidth]{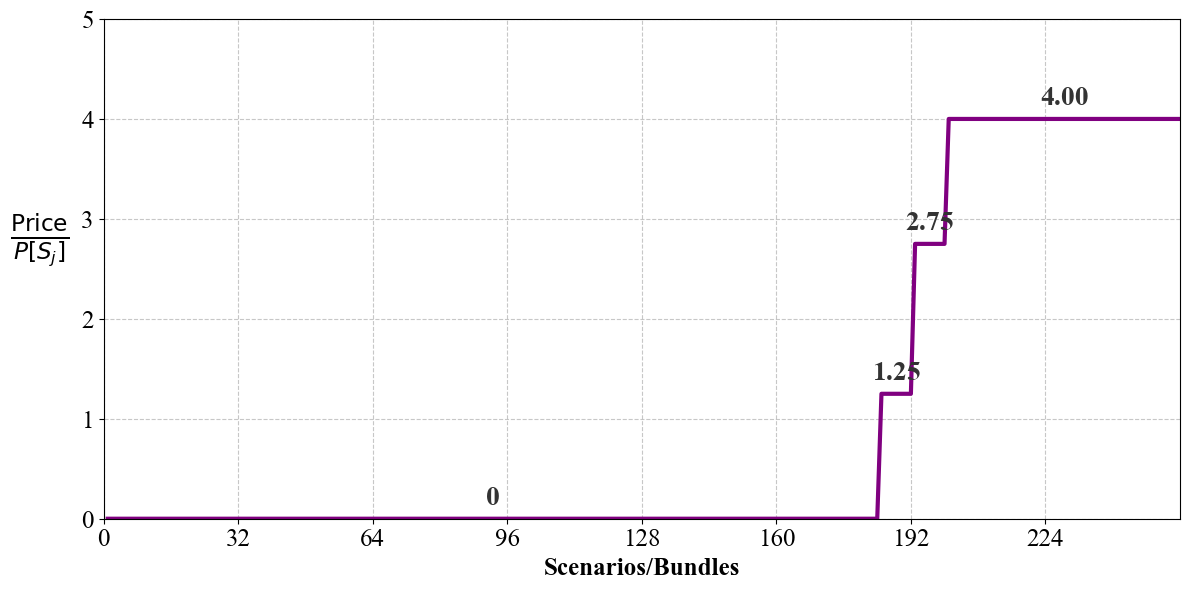}
        \caption{Bundle prices at iteration $m=5$.}
    \end{subfigure}
    \begin{subfigure}{0.30\textwidth}
        \centering
        \includegraphics[width=\linewidth]{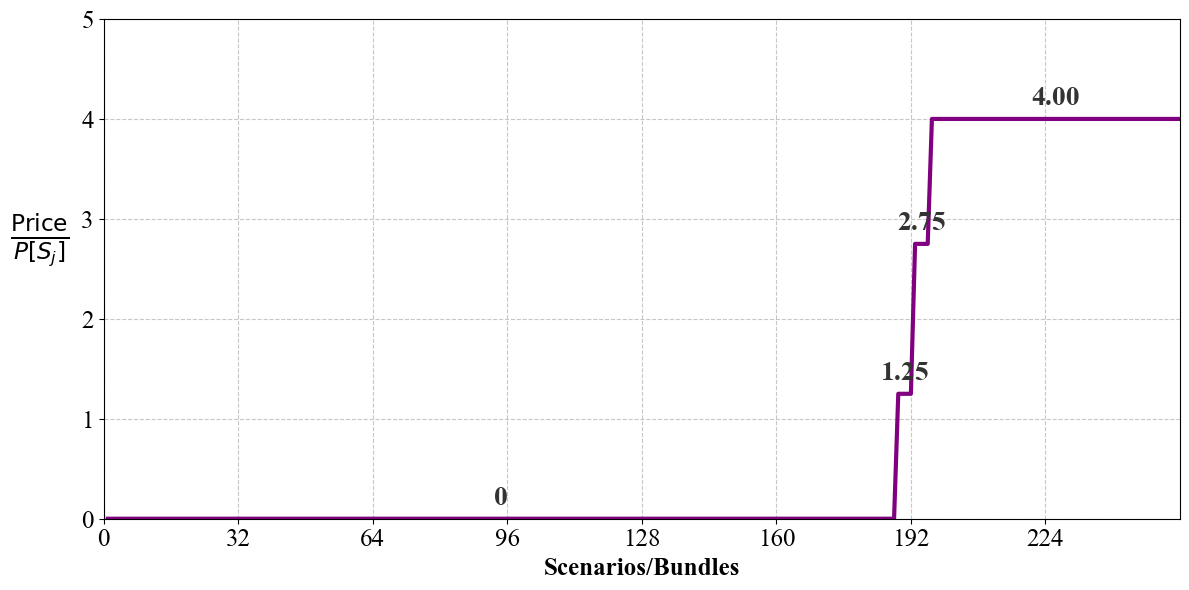}
        \caption{Bundle prices at iteration $m=6$.}
    \end{subfigure}
    
    \begin{subfigure}{0.3\textwidth}
        \centering
        \includegraphics[width=\linewidth]{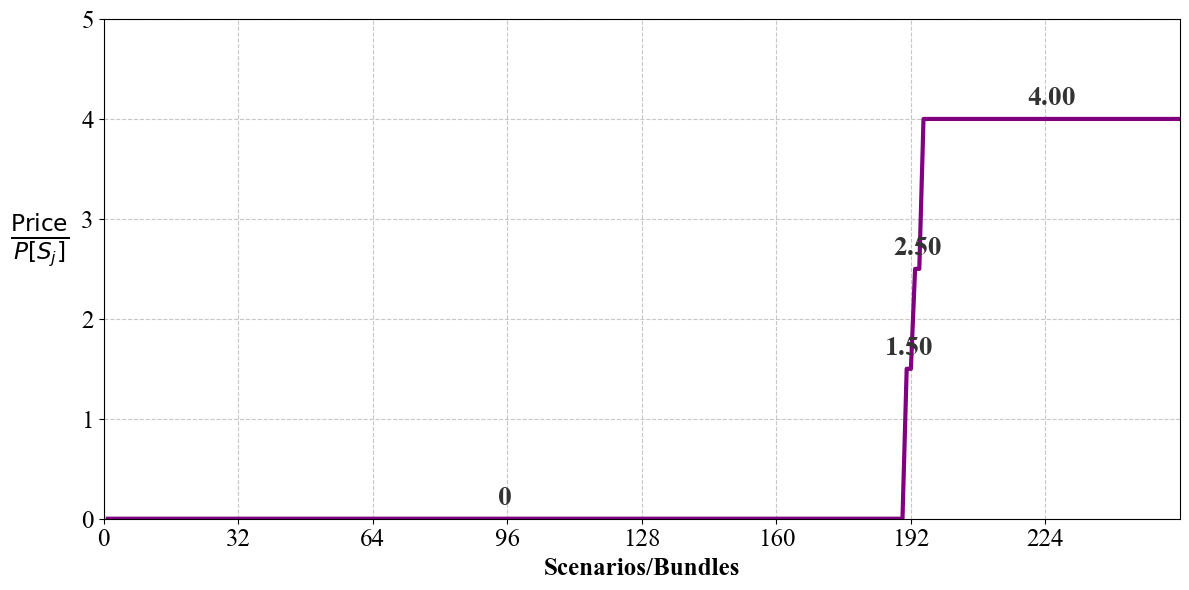}
        \caption{Bundle prices at iteration $m=7$.}
    \end{subfigure}
    \begin{subfigure}{0.3\textwidth}
        \centering
        \includegraphics[width=\linewidth]{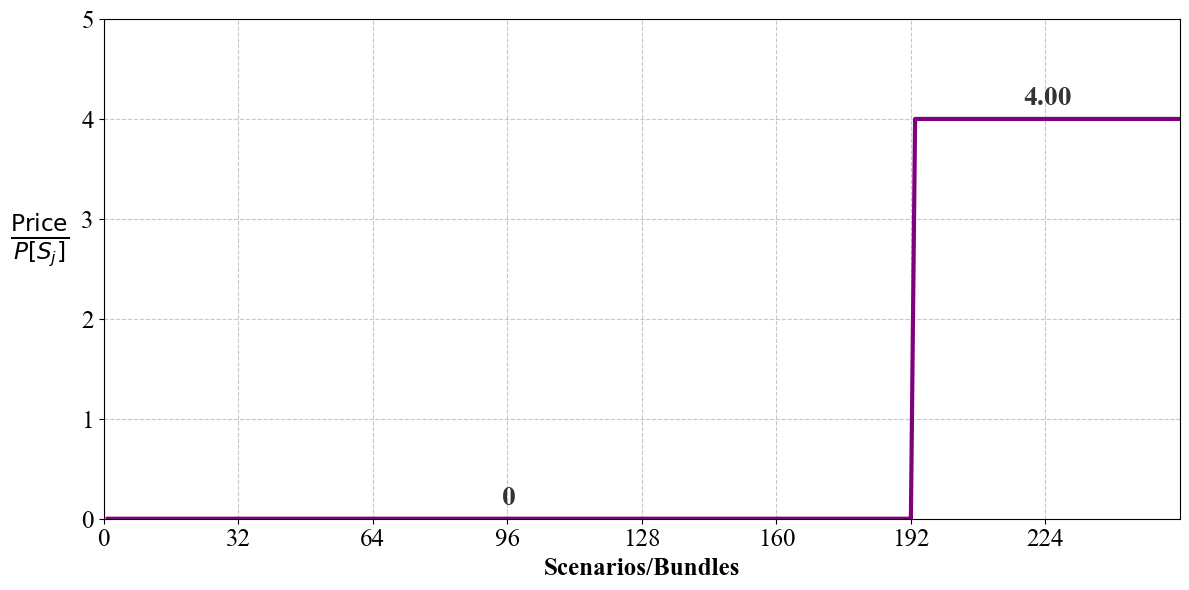}
        \caption{Risk price in a complete market.}
    \end{subfigure}
    
    \caption{Progression of risk price with bundle refinement.}
    \label{fig:five_images}
\end{figure}

One natural question is whether it is beneficial for the market to refine a bundle of scenarios with zero price. In other words, will scenarios within a zero-priced bundle eventually receive nonzero prices through successive refinements, or will their prices remain zero in all subsequent iterations? We provide a counterexample to illustrate that it is indeed possible for a scenario within a zero-priced bundle to obtain a strictly positive price after refinement. In Proposition \ref{prop:refinement_of_zero-priced_bundle}, we further characterize the circumstances under which a zero-priced bundle retains its zero price after refinement. 

\medskip
{\bf Example 2.}
Consider a competitive market with two risk-averse agents, each using $CVaR_{\beta_i}$ as their risk measure, with parameters $\beta_1 = 1/5$ and $\beta_2 = 1/2$ for agents 1 and 2, respectively. We assume 256 scenarios representing possible future outcomes, each with the same probability. Both agents are endowed with identical losses across scenarios, defined by: 
\[
z_{1s} = z_{2s} =
\begin{cases}
   \frac{s}{256}  & \text{for all } s\in\{1,...,256\} \text{ and } s \neq 64  \\
     \frac{s}{256} + 0.05  &  \text{for }s = 64
\end{cases}
\]
Assuming the same refinement method as described in Example 1, we can examine the price evolution of bundle containing scenario 64. Figure 4 illustrates how the price of such a bundle changes throughout the iterations of the refinement process. The bundle has a price of zero in iterations 1 and 2, and from iteration 3 onward, it begins to receive non-zero prices. 
\begin{figure}[htp!]
    \centering
    \includegraphics[width=0.65\textwidth]{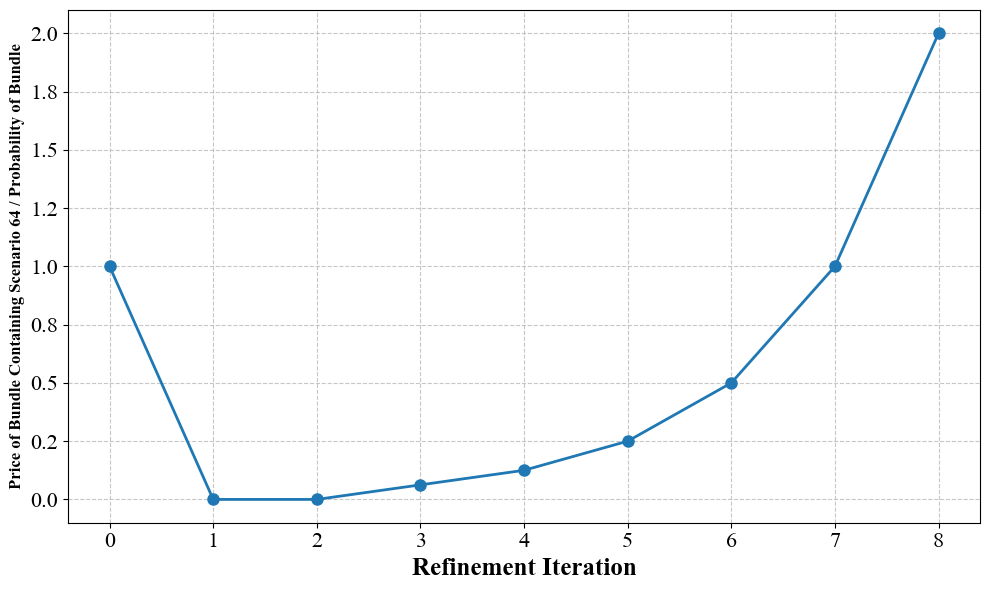}  
    \caption{Price progression of the bundle containing scenario 64 across refinement iterations.}
    \label{fig:figure4}
\end{figure}

However, it is important to note that the refinement process would maintain the price of the zero-priced bundle if the only bundle being refined is the one with a price of zero, while the other bundles remain unchanged. In other words, refining only the zero-priced bundle, while leaving the other bundles the same, would not result in an improvement in social welfare or provide any additional benefit to the market.
\begin{proposition}
\label{prop:refinement_of_zero-priced_bundle}
In a competitive market with $n$ risk-averse agents using \( CVaR_{\beta_i} \), refining only a zero-priced bundle  \( S_J \) into two new bundles \( S_{J_1} \) and \( S_{J_2} \), while leaving all other bundles unchanged, does not provide additional benefits to the market.
\end{proposition}
\begin{proof}
    Consider the primal formulation of the social optimization problem, where we use the linear representation of the $CVaR_{\beta_i}$ for each agent:
    \begin{align*}
        \min &\sum_{i=1}^n \{ \eta_i + \sum_{s=1}^{256}\frac{p_s}{\beta_i}u_{is}\}\\
         s.t. \  &\\
        & \eta_i + u_{is} \geq z_{is} - w_{is} && \forall i,s &[\pi_{is}]
        \\
        & w_{is} = t_{ij} && \forall i,\forall s\in S_j & [\alpha_{is}]\\
        & \sum_{i=1}^n t_{ij} = 0 && \forall j \quad & [\lambda_j]\\
         & u \geq 0
    \end{align*}
    where $\lambda_j$ is the price of bundle $j$ traded in the incomplete market. The dual problem for this linear program can be written as:
    \begin{align*}
        \max &\sum_{i=1}^n\sum_{s=1}^{256} z_{is} \pi_{is}\\
         s.t. \  &\\
        & \sum_{s=1}^{256} \pi_{is} = 1 && \forall i &[\eta_{i}]
        \\
        & \pi_{is} \leq \frac{p_s}{\beta_i} && \forall i,s & [u_{is}]\\
        & \sum_{s \in S_j} \pi_{is} = \lambda_j && \forall j \quad & [t_{ij}]\\
         & \pi \geq 0
    \end{align*}
 Now assume that in the incomplete market $\M_J$, bundle $S_J$  has a price of zero, i.e., \( \lambda_J = 0 \). We aim to show that the social welfare remains unchanged in the refined market $\M_{J_1 \cup J_2}$, where bundle $ S_J$ is split into two new bundles $S_{J_1}$ and $ S_{J_2}$, while the other bundles stay the same. To do this, observe that before the refinement, the following constraint holds for each agent:
 \begin{equation*}
       \sum_{s \in S_J} \pi_{is} = \lambda_J
  \end{equation*}
 After the refinement, the constraint becomes:
  \begin{equation*}
       \sum_{s \in S_{J_1}} \pi_{is} = \lambda_{J_1} \quad \text{and} \quad \sum_{s \in S_{J_2}} \pi_{is} = \lambda_{J_2}
  \end{equation*}
which makes the problem more constrained. Since $\lambda_J = 0$ in the market $\M_J$, the optimal solution satisfies $\pi_{is} = 0$ for all $i$ and all $s \in S_J$. This assignment also satisfies the new constraints with $\lambda_{J_1} = \lambda_{J_2} = 0$, and remains feasible in the refined market $M_{J_1 \cup J_2}$. Therefore, it is also optimal, implying that refining such a zero-priced bundle does not improve the market outcome.
\end{proof}

\section{Incomplete Markets of Bundled Instruments: discrete distributions with countably infinite states }
\label{Incomplete Markets of Bundled Instruments: discrete distributions with countably infinite states }

To extend the risk market beyond discrete distributions with finite scenarios, a natural approach is to consider cases where the scenarios are infinite but countable. Commonly used distributions in this context include the Poisson and geometric distributions. In this framework, the set of future states is $S = \mathbb{N}$, with $P$ representing the probability measure over it. Moreover, we will consider random variables that map $S$ to $\mathbb{R} $ and belong to $\L_1(S, P)$.

In this context, we will apply a variation of our refinement approach used in the finite case above, and we note that the process can continue infinitely, as in each iteration there exists a scenario bundle that can be further refined. The refinement yields (an infinite) sequence of incomplete markets $\{\M_n\}$, and the main question in this case is the asymptotic behavior of the sequence of incomplete markets $\{\M_n\}$. This problem can be tackled by examining whether the refinement process ensures convergence of the incomplete markets to the complete one in the limit. We can utilize Theorem 2 of \cite{Zame} by Zame to establish a convergence result of this nature. To this end, we employ the framework developed in \cite{Zame} and begin by reviewing some of its terminology and notation. In proving our convergence results, we rely on a particular approximation concept called \textit{uniform approximation} introduced in \cite{Zame}, where the uniform norm (sup norm) serves as the distance. For any $x, y \in \L_1(S, P)$, we can calculate the uniform distance of $x$ from $y$ with $d_\infty(x, y) = \sup_{s\in S}|x(s) - y(s)|$. This distance function defines a complete and metrizable topology on $\L_1(S, P)$, referred to as the \textit{uniform topology}. Thus, we can define the closure of a set $A \subset \L_1(S, P)$ with respect to the unifrom topology and denote it by $\text{cl}_\infty(A)$. For a given random variable $x$ and a set $A \subset \L_1(S, P)$, we say that $x$ can be \textit{uniformly approximated} by the elements of $A$ if and only if $x \in \text{cl}_\infty(A)$. 

Zame's converegence result is predicated on the concept of asymptotic completeness of the sequence of incomplete markets. Below, we remind the reader of the definition of this concept and state Theorem 2 of \cite{Zame}, adapted for our assumptions and for the reader's convenience.{\footnote{Zame also introduced another property, known as \textbf{asymptotic efficiency}, which is defined in terms of Pareto optimal allocations. However, this property does not provide additional insights in our case, as it is essentially equivalent to \textbf{asymptotic completeness}. This equivalence follows from the well-established result that Pareto optimal allocations are closely tied to welfare maximization, as outlined in the introduction and in particular, in our case where the agents' problems are convex optimization problems.}}

\begin{definition}
   The sequence $\{\M_n\}$ of the securities markets is said to be \underline{asymptotically complete} if, for every $\epsilon > 0$, there exists an integer $N_0$ such that for all $N > N_0$, every equilibrium allocation of $\M_n$ is within $\epsilon$ of an
equilibrium allocation of the underlying Walrasian economy $\M_{CM}$.
\end{definition}
\begin{theorem}
    Let $\{\alpha_n\} $ be the sequence of available bundles or securities in each incomplete market $\M_n$, with the first one being riskless, i.e., it pays off \$1 in all scenarios. The following are equivalent:
  \begin{enumerate}
      \item the sequence $\{\alpha_n\} $ is asymptotically complete;
      \item every Arrow-Debreu security can be uniformly approximated by the payoffs of a finite portfolio of security bundles $\alpha_n$.
  \end{enumerate}
\end{theorem}

Observe that any refinement method that induces asymptotic completeness of the markets will result in convergence of the incomplete market sequence to the complete market (using the above theorem). We will now lay out one such refinement scheme. We begin by defining the first bundle to encompass all possible scenarios. The second bundle is then defined to include all scenarios except the first, the third to exclude the first two scenarios, and this pattern continues accordingly. This approach clearly demonstrates that any given Arrow-Debreu security, such as the one paying \$1.00 in scenario $k$, can be replicated by investing \$1.00 in the $(k-1)^\text{th}$ instrument simultaneously shorting \$1.00 in the $k^\text{th}$ instrument. The payoff matrix for the refinement method would be as follows, where the element in row $i$ and column$j$ represents the payoff of asset $j$ in scenario $i$. 
\[
\begin{bmatrix}
1 & 0 & 0 & 0 & 0 & \cdots \\
1 & 1 & 0 & 0 & 0 & \cdots \\
1 & 1 & 1 & 0 & 0 & \cdots \\
1 & 1 & 1 & 1 & 0 & \cdots \\
\vdots & \vdots & \vdots & \vdots & \vdots & \ddots \\
\end{bmatrix}
\]
Market $\{\M_1\}$ is constructed by allowing participants to trade instrument 1, while $\{\M_2\}$ includes instruments $\{1,2\}$, and the market $\mathcal\{\M_j\}$ is constructed by including instruments $\{1, 2, \dots, j\}$ . We observed that it is guaranteed that this refining method can produce any Arrow-Debreu security in a finite number of iterations; thus, every Arrow-Debreu security can be uniformaly approximated by a finite portfolio of securities $\alpha_n$ available in the incomplete market $\M_n$. Therefore, the equilibrium of the sequence of incomplete markets $\{\M_n\}$ will consequently converge to the equilibrium of the complete market $\M_{CM}$.

\medskip
{\bf Example 3.} Consider a market with two risk-averse agents who evaluate their risk using $CVaR_{\beta_i}$, where $\beta_1 = 1-\exp(-1)$,  $\beta_2 = 1-2\exp(-1)$. The agents assume that the scenarios of future losses follow a Poisson distribution with a rate of occurrence $ \mu = 1 $ (implying a countably infinite set of possible scenarios indexed by $s \in \mathbb{N} \cup \{0\}$). Both agents are endowed with loss functions $Z_1$ and $Z_2$, where $Z_1: \mathbb{N} \cup \{0\} \to \mathbb{R}$ is defined as $Z_1(s) = 1$, and $Z_2: \mathbb{N} \cup \{0\} \to \mathbb{R}$ is defined as $Z_2(s) = s$. The market offers three financial instruments:
\begin{itemize}
    \item The \textbf{first} instrument pays \$1 in all scenarios.
    \item The \textbf{second} instrument pays \$1 in all scenarios except the first.
    \item The \textbf{third} instrument pays nothing in the first two scenarios and \$1 in all subsequent scenarios.
\end{itemize}
The equilibrium in the incomplete risk market is determined by solving the problem in (\ref{eq:Social Optimization}), which, under a linear formulation of $CVaR_\beta$, transforms into the following optimization problem:
    \begin{align*}
    \min \quad &\sum_{i=1}^2\{\eta_i + \sum_{s=0}^\infty \frac{P(S=s)U_i(s)}{\beta_i}\}\\ 
    s.t. \  &\\
    & \eta_i + U_i(0) \geq Z_i(0) - t_{i1} && \forall i
    \\ 
    & \eta_i + U_i(1) \geq Z_i(1) - t_{i1} - t_{i2} && \forall i\\
    & \eta_i + U_i(s) \geq Z_i(s) - t_{i1} - t_{i2}- t_{i3}  && \forall i,\forall s\geq2\\
& \sum_{i} t_{ij} = 0 &&\forall j\in \{1,2,3\}
    \end{align*}
Where \( t_{ij} \) represents the number of trade instruments of type \( j \) purchased by agent \( i \). Within the topological duality framework of Anderson \cite{Anderson1983}, the dual problem is given by:
\begin{align}
    \max \quad & \sum_{i=1}^2 \sum_{s=0}^\infty \pi_i(s) Z_i(s) = \sum_{s=0}^\infty \pi_1(s) + \sum_{s=0}^\infty s . \pi_2(s) \nonumber \\
    s.t. \  & \nonumber \\
    & \sum_{s=0}^\infty \pi_i(s) = 1  &&\forall i&\tag{i} \label{eq:constraint_i} \\
    &  \pi_i(s) \leq \frac{P(S = s)}{\beta_i}  &&\forall i ,\forall s\in\mathbb{N}\cup\{0\} &\tag{ii} \label{eq:constraint_ii} \\
    & \sum_{s=0}^\infty \pi_i(s) = \lambda_1  &&\forall i &\tag{iii} \label{eq:constraint_iii} \\ 
    & \sum_{s=1}^\infty \pi_i(s) = \lambda_2  &&\forall i &\tag{iv} \label{eq:constraint_iv} \\ 
    & \sum_{s=2}^\infty \pi_i(s) = \lambda_3  &&\forall i &\tag{v} \label{eq:constraint_v}
\end{align}
In this case,
\[
P(S = s) = \frac{\exp(-1)}{s!}
\]
Constraints (\ref{eq:constraint_i}) and (\ref{eq:constraint_ii}) represent agents' risk sets. Meanwhile, constraints (\ref{eq:constraint_iii}), (\ref{eq:constraint_iv}), and (\ref{eq:constraint_v}) not only determine instrument prices $\lambda_j$ but also ensure consensus among agents regarding the probabilities of scenarios 0 and 1. According to Theorems 1 and 2 of \cite{Anderson1983}, it can be shown that the primal feasible solutions,
\[
\eta = \begin{pmatrix} \eta_1 \\ \eta_2 \end{pmatrix} = \begin{pmatrix} 0 \\ 2 \end{pmatrix}, \quad
U_1 = \begin{pmatrix} 0 \\ 0 \\ 1 \\ 1 \\ 1 \\ 1 \\ \vdots \end{pmatrix}, \quad
U_2 = \begin{pmatrix} 0 \\ 0 \\ 0 \\ 1 \\ 2 \\ 3 \\ \vdots \end{pmatrix}, \quad
t_1 = -t_2 = \begin{pmatrix} t_{11} \\ t_{12} \\ t_{13} \end{pmatrix} = \begin{pmatrix} 1 \\ 0 \\ -1 \end{pmatrix}
\]
and dual feasible solutions,
\[
\pi_1(s) =
\begin{cases} 
0 & \text{if } s = 0 \\
\frac{\exp(-1)}{(1 - \exp(-1)) s!} & \text{if } s \geq 1
\end{cases} \quad \text{,}
\quad
\pi_2(s) =
\begin{cases} 
0 & \text{if } s = 0 \\
\frac{\exp(-1)}{1 - \exp(-1)} & \text{if } s = 1 \\
\frac{\exp(-1)}{(1 - 2\exp(-1)) 2} - \frac{\exp(-1)}{1 - \exp(-1)} & \text{if } s = 2 \\
\frac{\exp(-1)}{(1 - 2\exp(-1)) s!} & \text{if } s \geq 3
\end{cases}
\]
\[
\mathbf{\lambda} = \begin{pmatrix} \lambda_1 \\ \lambda_2 \\\lambda_3  \end{pmatrix} = \begin{pmatrix} 1 \\ 1 \\1 - \frac{\exp(-1)}{1-\exp(-1)}  \end{pmatrix}
\]
are also optimal for the primal and dual problems, respectively, due to the complementarity condition. Therefore the market prices of these instruments are approximately $\lambda \approx (1,1,0.418)$.

\section{Incomplete Markets of Bundled Trade Instruments: arbitrary probability space}
\label{Incomplete Markets of Bundled Trade Instruments: arbitrary probability space}
While a sufficiently large finite number of discrete scenarios is a an accepted
approximation to a continuum of outcomes, it is interesting to examine our successively refined approach to market completion in an arbitrary probability space. In this section, we will first lay out the necessary preliminaries to address our questions over an arbitrary probability space. Following an illustrative example, we will present our convergence results. 

\subsection{Preliminaries}

Our setting will be the space $\L_1(\Omega,\F,P)$ of integrable random variables (i.e. random variables $Z$ with $E[|Z|]<\infty$) on
an arbitrary probability space $(\Omega,\F,P)$.
This space contains the other spaces $\L_p(\Omega,\F,P)$ for $1\leq p\leq\infty$ (in nested order: $\L_{p_1}(\Omega,\F,P)\subseteq\L_{p_2}(\Omega,\F,P)$ if and only if $p_1\geq p_2$).
We will usually abbreviate $\L_p(\Omega,\F,P)$ to $\L_p$.
When dealing with random variables in $\L_1$ we will implicitly take
equality to be in the almost-sure sense; thus
we will usually write {\it e.g.} $Z_1=Z_2$ rather than $Z_1=Z_2$ a.s.

Suppose that $n\geq2$ agents each have an endowment (contingent liability),
with the endowment of agent $i$ ($i=1,\ldots,n$) being represented by the 
random variable $Z_i\in\L_p$ ($1\leq p<\infty$).The agents are able to trade risks: agent $i$ buys a
portfolio of securities with a payoff represented by the random variable
$W_i\in\L_p$, which has the effect of modifying his liability to $Z_i-W_i$.
Agents may trade only with each other: thus $\sumi W_i=0$.

We consider both complete and incomplete risk markets.
Let $(\F_m)_{m=1}^\infty$ be a filtration with
$\sigma\left(\bigcup_{m=1}^\infty \F_m\right)=\F$;
for each $m=1,2,\ldots$, the $\sigma$-field $\F_m$ represents an incomplete market
in which the only risks that may be traded are those with $W_i\in\F_m$.
The full $\sigma$-field $\F$ represents a complete market, in which
any risks may be traded.

\medskip
{\bf Example 4.} A useful example to bear in mind is $\Omega=[0,1]$ with
the Borel $\sigma$-field $\F$ and uniform (Lebesgue) probability measure $P$.
The $\sigma$-field $\F_m$ is that generated by partitioning $\Omega$ into
$2^m$ sub-intervals each of length $2^{-m}$. Similar to the finite, discrete case, this corresponds to an incomplete
risk market in which similar outcomes (here the nearby points in the unit
interval) are aggregated to make tradable bundles. There are $2^m$ distinct
bundles available to be traded; each agent attempts to acquire a portfolio
of them to ameliorate the risks with which they are endowed.

\medskip
{\bf Risk measures.}
The agents value their liabilities with coherent risk measures.
The risk measure used by agent $i$ is $r_i$, given by
\begin{equation}
 r_i(Z) = \hbox{max}\{E[ZY] : Y\in D_i\} \qquad\hbox{ for }Z\in\L_p
\label{eq:crm_definition}
\end{equation} 
where $D_i$ (the {\em risk set}) is a closed convex subset of
$D_0=\{Y\in\L_1: Y\geq0\hbox{ and }E[Y]=1\}$.
The elements of $D_0$ may be thought of as probability density functions,
or Radon-Nikodym derivatives with respect to the base measure $P$.

\medskip
{\bf Assumptions.} We make the following assumptions throughout.
\begin{enumerate}
\item
Each $D_i$ is a closed and bounded set in $\L_q$
(where $p^{-1}+q^{-1}=1$). That is, $\sup\{\|Y\|_q: Y\in D_i\}<\infty$.
Equivalently, the risk measure $r_i$ takes finite values for any
risks in $\L_p$.
\item 
$\bigcap_{i=1}^n D_i$ is non-empty.
\end{enumerate}

\medskip
{\bf Interpretation.} The optimum $Y^*$ in (\ref{eq:crm_definition})
can be interpreted as a ``price of risk'', in the sense that agent $i$
would be willing to pay $E[WY]$ to reduce his liability from $Z$ to $Z-W$.

\medskip
{\bf Example 5.} For the $\hbox{CVaR}_\beta$ risk measure ($0<\beta<1$),
the risk set is $\{Y\in D_0:0\leq Y\leq \beta^{-1}\}$.
If all agents are using CVaR risk measures (possibly with different values
of $\beta$), we can take $p=1$ and $q=\infty$.

\medskip
{\bf Example 6.} For the Good-Deal risk measure,  see, e.g., \cite{Cochrane-Saa-Requejo-2000}, with base measure $P$
and parameter $\nu>1$,
the risk set is $\{Y\in D_0:E[Y^2]\leq\nu^2\}$.
If all agents are using this type of risk measure (possibly with
different values of $\nu$), we can take $p=q=2$.

\medskip
{\bf The weak* topology.}
It will sometimes be convenient to equip $D_i$ with the weak* topology
inherited from $\L_q$ in its role as the dual space of $\L_p$.
Since $D_i$ is  closed and convex, it is closed in the weak* topology as well as in the
$\L_q$ and $\L_1$ norm topologies. Since $D_i$ is also bounded in $\L_q$,
it is weak* compact by the Banach-Alaoglu theorem.
Among other things, this assures us that the weak*-continuous function
$Y\mapsto E[XY]$ in (\ref{eq:crm_definition}) attains its maximum over $D_i$
at some $Y^*\in D_i$.

\subsection{The social optimization problem}

The market clearing problem can be expressed as a social optimization:
\begin{equation}
\hbox{min}\quad
\sumi r_i(Z_i-W_i)
 \quad\hbox{s.t.}\quad
W_1,\ldots,W_n\in\L_p(\F_m)\hbox{ and }\sumi W_i=0 .
\label{eq:social_opt}
\end{equation} 
For the complete market, the space $\L_p(\F_m)$ is replaced by $\L_p$.

\medskip
{\bf Remark.} Since $\F_m\subseteq\F_{m+1}$, Problem (\ref{eq:social_opt})
becomes less constrained as $m$ increases. Its optimal value---the optimal
social welfare---thus increases with $m$. We can regard (\ref{eq:social_opt})
as representing a succession of incomplete markets, each of which improves on
the last by adding further instruments of trade in an attempt to approximate
the complete market represented by $\F$.

\medskip
Making use of (\ref{eq:crm_definition}) expands (\ref{eq:social_opt}) to
$$ \hbox{min}\quad
\sumi \hbox{max}\left\{E[(Z_i-W_i)Y]: Y\in D_i\right\}
 \quad\hbox{s.t.}\quad
W_1,\ldots,W_n\in\L_p(\F_m)\hbox{ and }\sumi W_i=0 ,
  $$
or equivalently
\begin{equation}
 \hbox{min}\quad
 \hbox{max}\left\{\sumi E[(Z_i-W_i)Y_i]: Y_i\in D_i\ \forall i\right\}
 \quad\hbox{s.t.}\quad
W_1,\ldots,W_n\in\L_p(\F_m)\hbox{ and }\sumi W_i=0 .
\label{eq:min_max}
\end{equation}

\begin{theorem}
For any elements $W_1^*,\ldots,W_n^*$ optimal for the outer problem of (\ref{eq:min_max}), elements $Y_i\in D_i$ optimal for the corresponding inner
problem are given by any optimum of the problem
\begin{equation}
 \hbox{max}\quad
 \sumi E[Z_i Y_i]
 \quad\hbox{s.t.}\quad
 Y_i\in D_i \ \forall i \hbox{ and $E[Y_i|\F_m]$ is the same for all $i$.} 
\label{eq:dual_Fm}
\end{equation}
\label{th:zy_thm}
\end{theorem}
In the complete-market case, (\ref{eq:dual_Fm}) requires the $Y_i$ themselves
to be the same for all $i$, and so reduces to
\begin{equation}
 \hbox{max}\quad
  E\left[\left(\sumi Z_i\right) Y_0\right]
 \quad\hbox{s.t.}\quad
 Y_0\in \bigcap_{i=1}^n D_i.
\label{eq:dual_complete}
\end{equation}
Problem (\ref{eq:dual_complete}) exhibits the consolidation effect
discussed in \cite{Ralph-Smeers-2015}: since a complete market can re-assign liabilities
between agents arbitrarily, it does not matter which agent has which
liability initially, but only what the total liability is.

\medskip
{\bf Interpretation.} The optimum $Y_i$ in (\ref{eq:dual_Fm}) also has
an interpretation as a price of risk:
after the inter-agent risk market has cleared, agent $i$ would now be
willing to pay $E[WY_i]$ to reduce his liability from $Z_i-W_i$ to
$Z_i-W_i-W$. 
The condition that $E[Y_i|\F_m]$ is the same for all $i$ is equivalent to
$E[WY_i] = E[WY_j]$ for any two agents $i,j$ and any $W\in\F_m$.
This can be interpreted as a requirement that after trading ceases,
agents must agree on the price of any tradable risk
(for if not, further trades would be possible).
Non-tradable risks may still be priced differently by different agents, however.
In the complete-market case, we have the familiar result that a single
pricing scheme, shared by all agents and applicable to all risks, emerges.

\begin{proof}[Proof of Theorem \ref{th:zy_thm}.]
To simplify notation, we write (\ref{eq:min_max}) more compactly as
\begin{equation}
 \hbox{min}\quad
 \hbox{max}\left\{E[(Z-W)\cdot Y]: Y\in D\right\}
 \quad\hbox{s.t.}\quad W\in S ,
\label{eq:min_max_dot}
\end{equation}
where $W=(W_1,\ldots,W_n)\in\L_p^n$ (and similarly $Z$),
$Y=(Y_1,\ldots,Y_n)\in\L_q^n$,
$\cdot$ denotes the usual dot product,
$D=\prod_{i=1}^n D_i$,
and $S$ is the linear subspace of $\L_p^n$ given by
$S=\left\{(W_1,\ldots,W_n)\in\L_p(\F_m)^n : \sumi W_i=0\right\}$.

Let $v\in{\mathbb R}$ be the optimal value of (\ref{eq:min_max_dot});
the assumed existence of an optimum implies that $v$ is a finite value.
Equip each $D_i$ with the weak* topology inherited from $\L_q$ in its role as
the dual space of $\L_p$. Equip $D$ with the product of these topologies,
making $D$ a compact Hausdorff space.
Then by Theorem 2 of \cite{Fan}, 
\begin{equation*}
\begin{split}
v &= \inf\left\{\hbox{max}\left\{E[(Z-W)\cdot Y]: Y\in D\right\}: W\in S\right\}\\
  &= \hbox{max}\left\{\inf\left\{E[(Z-W)\cdot Y]: W\in S\right\}: Y\in D\right\} .
\end{split}
\end{equation*}
Let $S^\perp = \{Y\in\L_q^n:E[W\cdot Y]=0\ \forall W\in S\}$.
For any $Y\notin S^\perp$, we have
$\inf\left\{E[(Z-W)\cdot Y]: W\in S\right\}=-\infty$ since $S$ is a linear subspace.
Hence
\begin{equation*}
v = \hbox{max}\left\{\inf\left\{E[(Z-W)\cdot Y]: W\in S\right\}:
         Y\in D\cap S^\perp\right\}
\end{equation*}
or, since $E[W\cdot Y]=0$ when $W\in S$ and $Y\in S^\perp$,
\begin{equation}
\begin{split}
v &= \hbox{max}\left\{\inf\left\{E[Z\cdot Y]: W\in S\right\}:
         Y\in D\cap S^\perp\right\}\\
  &= \hbox{max}\left\{E[Z\cdot Y]: Y\in D\cap S^\perp\right\} .
\end{split}
\label{eq:dual_Fm_dot}
\end{equation}
Let $Y^*$ be any element of $D\cap S^\perp$ attaining this maximum value.
(The existence of $Y^*$ is assured because $S^\perp$ is weak*-closed, and so
$D\cap S^\perp$ is weak*-compact.)
Let $W^*$ be optimal for (the outer problem of) (\ref{eq:min_max_dot}).
Then $E[(Z-W^*)\cdot Y^*] = E[Z\cdot Y^*] = v$.
That is, $Y^*$ achieves the optimal value $v$ in the problem
$$  \hbox{max}\quad
 E[(Z-W^*)\cdot Y] \quad\hbox{s.t.}\quad Y\in D ,
  $$
which is the inner optimization of (\ref{eq:min_max_dot}).

It remains to show that (\ref{eq:dual_Fm_dot}) is the same problem as
(\ref{eq:dual_Fm}), {\it i.e.} that
\begin{equation*}
S^\perp = \left\{(Y_1,\ldots,Y_n)\in\L_q^n: E[Y_i|\F_m]\hbox{ is the same for all }i\right\} .
\end{equation*}
To see this, observe that if $W\in S$ and $Y\in S^\perp$,
\begin{equation}
\begin{split}
0 &= E[W\cdot Y] = E\left[\sumi W_i Y_i\right]
= E\left[\sumi E\left[W_i Y_i|\F_m\right]\right]
= E\left[\sumi W_i E\left[Y_i|\F_m\right]\right]\\
&= E\left[\sum_{i=1}^{n-1} W_i E\left[Y_i|\F_m\right]
   + \left(-\sum_{i=1}^{n-1}W_i\right)E\left[Y_n|\F_m\right]\right]\\
&= E\left[\sum_{i=1}^{n-1} W_i \left(E\left[Y_i|\F_m\right] - 
                                    E\left[Y_n|\F_m\right]\right)\right] .
\end{split}
\label{eq:wy_perp}
\end{equation}
Hence, for any $i=1,\ldots,n-1$ and any $U\in\L_p(\F_m)$,
$$ E\left[U\left(E\left[Y_i|\F_m\right] - E\left[Y_n|\F_m\right]\right)\right] = 0
  $$
(by taking $W_i=U$, $W_n=-U$).
It follows that $E\left[Y_i|\F_m\right] = E\left[Y_n|\F_m\right]$.
Conversely, if $E\left[Y_i|\F_m\right]$ is the same for all $i$, then
(\ref{eq:wy_perp}) demonstrates that $Y\in S^\perp$.
\end{proof}

\medskip
{\bf Example 7.}
Consider a market with two agents who are endowed with losses $Z_1$ and $Z_2$, which are random variables in $\L_1(\Omega,\F,P)$. The sample space is defined as $\Omega = [0,1]$, where $\F $ is the Borel $\sigma$-algebra and the probability measure $P$ is uniform. The losses are specified as $Z_1: \omega \mapsto 2\omega$ and $Z_2: \omega \mapsto 1-\omega$. Additionally, $\F_1 = \{\Omega, [0,\frac{1}{2}), [\frac{1}{2},1], \emptyset\}$ represents an incomplete market in which only two products are
traded: the bundles $[0,\frac12)$ and $[\frac12,1]$, representing ``smaller''
and ``larger'' outcomes respectively. 
The agents evaluate their risk with $CVaR_{\beta_i}$.
The risk market prices
can be found by solving the problem (\ref{eq:dual_Fm}), which here takes the form

\begin{equation}
\begin{split}
\max \quad & \sum_{i=1}^2 E[Z_i Y_i] \\
\text{s.t.} \quad & 0 \leq Y_i \leq \beta_i^{-1}, \quad E[Y_i] = 1 \quad \text{for} \ i = 1,2, \\
& E[Y_1|\F_1] = E[Y_2|\F_1].
\end{split}
\label{eq:two_agent_cts_optimization}
\end{equation}

Note that $Z_1$ is increasing and $Z_2$ decreasing on $[0,1]$.
Hence, we need consider only $Y_1$, $Y_2$ of the form
\begin{equation}
\begin{split}
Y_1 &= \beta_1^{-1} \left( 1_{[a, \frac{1}{2})} + 1_{[b, 1]} \right) \\
\text{and} \quad
Y_2 &= \beta_2^{-1} \left( 1_{[0, c)} + 1_{[\frac{1}{2}, d]} \right) .
\end{split}
\label{eq:two_agent_abcd}
\end{equation}

since for any other probability densities $Y_1$, $Y_2$ the objective of
(\ref{eq:two_agent_cts_optimization}) could be improved by rearranging the
probability mass on $[0,\frac12)$ and/or $[\frac12,1]$.
The condition $E[Y_1]=E[Y_2]=1$ gives us
$\beta_1^{-1}(\frac12 - a + 1 - b) = \beta_2^{-1}(c+d-\frac12)=1$,
while the condition $E[Y_1|\F_1]=E[Y_2|\F_1]$ gives us
$\beta_1^{-1}(\frac12 - a)=\beta_2^{-1} c$ and
$\beta_1^{-1}(1-b)=\beta_2^{-1}(d-\frac12)$.
With these conditions, (\ref{eq:two_agent_abcd}) reduces to the
one-parameter family

\begin{equation}
\begin{split}
Y_1 &= \beta_1^{-1} \left( 1_{[\frac{1}{2} - \beta_1 u, \frac{1}{2})} + 1_{[1 - \beta_1(1 - u), 1]} \right) \\
\text{and} \quad
Y_2 &= \beta_2^{-1} \left( 1_{[0, \beta_2 u)} + 1_{[\frac{1}{2}, \frac{1}{2} + \beta_2(1 - u)]} \right) .
\end{split}
\label{eq:two_agent_split}
\end{equation}

where $u\geq0$.
We must have
$$ \beta_1 u \leq\frac12, \quad
   \beta_1 (1-u) \leq\frac12, \quad
   \beta_2 u \leq\frac12, \quad\hbox{ and }\quad
   \beta_2 (1-u) \leq\frac12 ,
  $$
that is,
$$ \hbox{max}(u,1-u) \leq \frac{1}{2\;\hbox{max}(\beta_1,\beta_2)} ,
  $$
or
$$ 1-\lambda \;\leq\; u \;\leq\; \lambda,
 \qquad\hbox{ where }\quad
 \lambda = \frac{1}{2\;\hbox{max}(\beta_1,\beta_2)} .
  $$
Now consider the objective of (\ref{eq:two_agent_cts_optimization}) :
\begin{equation*}
\begin{split}
v &= \sum_{i=1}^2 E[Z_i Y_i] \\
  &= \beta_1^{-1} \int_{\frac12-\beta_1 u}^{\frac12} 2x \, dx + \beta_1^{-1} \int_{1-\beta_1(1-u)}^{1} 2x \, dx \\
  &\quad + \beta_2^{-1} \int_{0}^{\beta_2 u} (1-x) \, dx + \beta_2^{-1} \int_{\frac12}^{\frac12 + \beta_2(1-u)} (1-x) \, dx.
\end{split}
\end{equation*}

Differentiating yields
\begin{equation}
\label{eq:dvdu}
\frac{dv}{du} \;=\;
-\frac12 + (2\beta_1 + \beta_2)(1-2u) .
\end{equation}

We distinguish two cases of interest.
Firstly, if $2\beta_1+\beta_2\leq\frac12$, then $\lambda>1$ and
$\frac{dv}{du}\big|_{u=0}\leq0$. In this case the maximum occurs for $u=0$,
meaning that the optimal probability densities $Y_1^*$ and $Y_2^*$ both
place all of their probability mass on $[\frac12,1]$. The other bundle
$[0,\frac12)$ has zero market value.

Secondly, if $2\beta_1+\beta_2>\frac12$ but
$\hbox{max}(\beta_1,\beta_2)\leq\frac12$, then $\lambda\geq1$ and
$\frac{dv}{du}\big|_{u=0}\geq0$. So, the maximum can be found by
setting the derivative in (\ref{eq:dvdu}) to zero, giving
$$ u = \frac12 - \frac{1}{4(2\beta_1+\beta_2)} .
  $$
In this case both bundles have positive market prices, given by
$$ E[Y_1^*|\F_1] = E[Y_2^*|\F_1]
 = \begin{cases}
 1 - \frac{1}{2(2\beta_1+\beta_2)} &, \hbox{ for } [0,\frac12)\\
 1 + \frac{1}{2(2\beta_1+\beta_2)} &, \hbox{ for } [\frac12,1] .
 \end{cases}
  $$

\subsection{Convergence of risk prices}

In this section we show that the price of risk discovered by a complete
market can be approximated by the prices discovered on incomplete
markets, in the sense that the incomplete-market prices $Y_i$ on a market
$\F_m$ converge to the complete-market price $Y_0$ as $m\to\infty$.

Problem (\ref{eq:dual_Fm}) is a relaxation of Problem (\ref{eq:dual_complete}).
The degree of relaxation depends on the parameter $m$: the larger the
value of $m$, the slighter the relaxation. This suggests that
(\ref{eq:dual_Fm}) should also be regarded as a perturbation of
(\ref{eq:dual_complete}); in this light, our aim is to show that the set of
optimal solutions---the price of risk---behaves continuously with respect to
such a perturbation.

More precisely, the concept required to describe this behaviour is
{\em upper semicontinuity} of the optimal solution set.

\begin{definition}
Let $X$ and $M$ be Hausdorff topological spaces.
A multifunction $S:M\to2^X$ is {\em upper semicontinuous} at a point
$m_0\in M$ if for every neighbourhood $V$ of the set $S(m_0)$ there is a
neighbourhood $U$ of $m_0$ such that $m\in U\implies S(m)\subseteq V$.
\end{definition}

Intuitively, upper semicontinuity means that the perturbation from $m_0$
to $m$ does not introduce elements of $S(m)$ which are very different from
any element of $S(m_0)$.

The main result of this section can now be stated as follows.

\begin{theorem}
\label{th:prices_are_usc}
Equip $M={\mathbb N}\cup\{\infty\}$ with the usual topology (the
neighbourhoods of $\infty$ are those sets whose complements are finite
subsets of $\mathbb N$).
Equip $\L_q^n$ ($1<q\leq\infty$) with its weak* topology as the dual space
of $\L_p^n$.
Let $S(m)$ be the set of optimal solutions
of (\ref{eq:dual_Fm}) if $m<\infty$, or of (\ref{eq:dual_complete}) if
$m=\infty$.
Then the multifunction $S$ is upper semicontinuous at $\infty$.
\end{theorem}

We defer the proof to an appendix.
The result of Theorem \ref{th:prices_are_usc} can be stated in less abstract
language as follows.

\begin{corollary}
\label{cor:prices_converge}
Suppose that $Y_0$ is an optimal solution of (\ref{eq:dual_complete}),
and also that for each $m<\infty$,
$(Y_1^m,\ldots,Y_n^m)$ is an optimal solution of (\ref{eq:dual_Fm}).
Then for any $W_1,\ldots,W_n\in\L_p$ with $\sumi W_i=0$ we have
$$ \sumi E[W_i Y_i^m] \;\to\; \sumi E[W_i Y_0] \;=\; 0
  \quad\hbox{ as }m\to\infty .
  $$
\end{corollary}

\begin{proof}[Proof of Corollary \ref{cor:prices_converge}.]
For any $\epsilon>0$, let
$V=\left\{(Y_1,\ldots,Y_n):\left|\sumi E\left[W_i(Y_i - Y_0)\right]\right|<\epsilon\right\}$, a neighbourhood of $S(\infty)$ in the weak* topology.
By upper semicontinuity, there exists $m_1\in{\mathbb N}$ such that
for $m\geq m_1$ we have $S(m)\subseteq V$, and in particular
$(Y_1^m,\ldots,Y_n^m)\in V$, so that
$\left|\sumi E\left[W_i Y_i^m\right]\right|<\epsilon$.
\end{proof}

\appendix 
\section{Proof of Theorem \ref{th:prices_are_usc}.}
Theorem \ref{th:prices_are_usc} can be derived from the following
result, which appears as Proposition 4.4 in \cite{BS}.

\begin{proposition}
\label{prop:BS44}
Let $X$ and $M$ be Hausdorff topological spaces, and let $m_0\in M$.
Let the multifunctions $\Phi:M\to2^X$ and $S:M\to2^X$ give
(respectively) the feasible and the optimal sets of the parameterized
optimization problem
\begin{equation}
 \hbox{max}\quad f(x,m) \quad\hbox{s.t.}\quad x\in\Phi(m) .
\end{equation}
Suppose that
\begin{enumerate}
\item the function $f$ is continuous on $X\times U$;
\item for any neighbourhood $V$ of the set $S(m_0)$ there exists a
neighbourhood $U$ of $m_0$ such that
$V\cap\Phi(m)\neq\emptyset\ \forall m\in U$;
\item $\exists\alpha\in{\mathbb R}$ and a compact set $C\subseteq X$
such that for every $m$ in a neighbourhood of $m_0$, the level set
$$ \{x\in\Phi(u): f(x,m)\geq\alpha\}
  $$
is a non-empty subset of $C$;
\item $\{(x,m)\in X\times U: x\in\Phi(m)\}$ is closed.
\end{enumerate}
Then
\begin{enumerate}
\item the optimal value function $\nu(m)$ is continuous at $m_0$, and
\item $S$ is upper semicontinuous at $m_0$.
\end{enumerate}
\end{proposition}

\begin{proof}[Proof of Theorem \ref{th:prices_are_usc}.]
We verify the sufficient conditions given in Proposition \ref{prop:BS44}.
The first three conditions are straightforward.
\begin{enumerate}
\item 
The objective function of (\ref{eq:dual_Fm}) and (\ref{eq:dual_complete})
is not perturbed: it is $(Y_1,\ldots,Y_n)\mapsto \sumi E[Z_i Y_i]$
for all $m$. (In (\ref{eq:dual_complete}), the $Y_i$ are constrained to
be all equal.) As the objective is a linear functional on $\L_q^n$, it
is continuous with respect to the weak* topology.
\item
This condition is easily satisfied: since (\ref{eq:dual_Fm}) is a
relaxation of (\ref{eq:dual_complete}), we have $S(\infty)\subseteq\Phi(m)$
for any $m<\infty$.
\item 
Let $C=D_1\times\cdots\times D_n$. Then $C$ is a closed, bounded, and convex
subset of $\L_q^n$. A convex set is strongly closed iff it is weakly closed,
and so $C$ is closed in the weak* topology also. It follows by the
Banach-Alaoglu theorem that $C$ is compact in the weak* topology.
Let $Y_0$ be any optimal solution of (\ref{eq:dual_complete}); such a
solution exists as we have assumed $\bigcap_{i=1}^n D_i$ is non-empty.
Take $\alpha = E\left[\left(\sumi Z_i\right) Y_0\right]$.
Then the level set
$$ \left\{(Y_1,\ldots,Y_n): Y_i\in D_i\ \forall i,
 E[Y_i|\F_m] \hbox{ are equal }\forall i,
 \hbox{ and } \sumi E[Z_i Y_i]\geq\alpha\right\}
  $$
is non-empty (since it includes $(Y_0,\ldots,Y_0)$) and contained in $C$.
\end{enumerate}

For the fourth condition, it is required to show that the set
$$ K = \left\{(Y_1,\ldots,Y_n,m): Y_i\in D_i\ \forall i,
  E[Y_i|\F_m] \hbox{ are equal }\forall i,
 \hbox{ and } m\in M
      \right\}
  $$
is closed in $\L_q^n\times M$.
Here we are using the notation $\F_\infty$ for the full $\sigma$-field $\F$
to include the complete-market case $m=\infty$.

Suppose that the sequence $(Y_1^k,\ldots,Y_n^k,m_k)$ in $K$ converges to a
limit $(Y_1^\infty,\ldots,Y_n^\infty,m_\infty)$ in $\L_q^n\times M$
as $k\to\infty$.
It will suffice to show that $(Y_1^\infty,\ldots,Y_n^\infty,m_\infty)\in K$.
As noted above, $D_1\times\cdots\times D_n$ is weak*-closed in $\L_q^n$,
and so we immediately have $Y_i^\infty\in D_i\ \forall i$.

The convergence is in the weak* sense: for any $W_1,\ldots,W_n\in\L_p$,
$\sumi E[(Y_i^k - Y_i^\infty)W_i]\to0$ as $k\to\infty$.
If we further suppose that $\sumi W_i=0$, then we have
\begin{equation}
\label{eq:m_finite_and_infinite}
\begin{split}
\sumi E[Y_i^k W_i] \;
&=\; \sumi E\left[E\left[Y_i^k(W_i - E[W_i|\F_{m_k}] + 
                         E[W_i|\F_{m_k}])|\F_{m_k}\right]\right]\\
&=\; \sumi E\left[E\left[Y_i^k(W_i - E[W_i|\F_{m_k}])|\F_{m_k}\right]\right]+\\
&\; \sumi E\left[E\left[Y_i^k|\F_{m_k}] E[W_i|\F_{m_k}\right]\right]\\
&=\; \sumi E\left[Y_i^k(W_i - E[W_i|\F_{m_k}])\right]
 + E\left[E[Y_1^k|\F_{m_k}] E\left[\sumi W_i\Bigg|\F_{m_k}\right]\right]\\
&=\; \sumi E\left[Y_i^k(W_i - E[W_i|\F_{m_k}])\right] ,
\end{split}
\end{equation}
since $E[Y_i^k|\F_{m_k}]=E[Y_1^k|\F_{m_k}]\ \forall i$.

We now distinguish the cases $m_\infty=\infty$ and $m_\infty<\infty$.
In the case where $m_\infty=\infty$ (and so $m_k\to\infty$), we can observe
that $E[W_i|\F_{m_k}]\to W_i$ in $\L_p$ as $k\to\infty$.
For $p=1$ this is the L\'evy martingale convergence theorem;
for $p>1$ we rely instead on the martingale $\L_p$ convergence theorem
(\cite{Durrett}; the L\'evy martingale $(E[W_i|\F_{m_k}])_{k=1}^\infty$ is bounded in $\L_p$, since $E[|E[W_i|\F_{m_k}]|^p]\leq E[|W_i|^p]$ by Jensen's inequality).
Applying the H\"older inequality to (\ref{eq:m_finite_and_infinite}) gives
$$ \left|\sumi E[Y_i^k W_i]\right|
 \leq \sumi \|Y_i^k\|_q \|W_i - E[W_i|\F_{m_k}]\|_p
 \leq c \sumi \|W_i - E[W_i|\F_{m_k}]\|_p
 \to 0 \hbox { as }k \to\infty ,
  $$
where $c = \sup\left\{\|Y\|_q: Y\in\bigcup_{i=1}^n D_i\right\} < \infty$.
As we have $\sumi E[Y_i^k W_i] \to \sumi E[Y_i^\infty W_i]$, we may
conclude that $\sumi E[Y_i^\infty W_i] = 0$.
In particular, for any indices $i,j$ with $i\neq j$ and any $W\in\L_p$,
we can let $W_i=W$, $W_j=-W$, and $W_\ell=0$ for $\ell\notin\{i,j\}$
to obtain $E[(Y_i^\infty - Y_j^\infty)W]=0$, from which it follows that
$Y_i^\infty = Y_j^\infty$ in $\L_q$; thus the sequence limit
$(Y_1^\infty,\ldots,Y_n^\infty,\infty)$ belongs to $K$.

Turning to the case $m_\infty < \infty$, we note that in this case
$m_k=m_\infty$ for all but finitely many $k$; we can thus simplify the notation
by assuming, without loss of generality, that $m_k=m_\infty=m\ \forall k$.
Equation (\ref{eq:m_finite_and_infinite}) then gives us that for any
$W_1,\ldots,W_n\in\L_p$ with $\sumi W_i=0$ and $W_i\in\F_m$
(so that $E[W_i|\F_m]=W_i$) we have $\sumi E[Y_i^k W_i]=0$.
Since $\sumi E[Y_i^k W_i] \to \sumi E[Y_i^\infty W_i]$, it follows that
$\sumi E[Y_i^\infty W_i]=0$. Hence
\begin{equation}
\label{eq:m_finite}
 \sumi E\left[E\left[Y_i^\infty|\F_m\right]W_i\right]
 = \sumi E\left[E\left[Y_i^\infty W_i|\F_m\right]\right]
 = \sumi E\left[Y_i^\infty W_i\right]
 = 0 .
\end{equation} 
But this last conclusion also holds without the requirement that $W_i\in\F_m$,
since in general
\begin{equation}
\begin{split}
\sumi E\left[E\left[Y_i^\infty|\F_m\right]W_i\right] \;
&=\; \sumi E\left[E\left[Y_i^\infty|\F_m\right](W_i
   - E[W_i|\F_m] + E[W_i|\F_m])\right]\\
&=\; \sumi E\left[E\left[E\left[Y_i^\infty|\F_m\right](W_i - E[W_i|\F_m])|\F_m\right]\right]+\\
&\; \sumi E\left[E\left[Y_i^\infty|\F_m\right]E[W_i|\F_m]\right]\\
&=\; \sumi E\left[E\left[Y_i^\infty|\F_m\right]E\left[W_i - E[W_i|\F_m]|\F_m\right]\right]+\\
 &\; \sumi E\left[E\left[Y_i^\infty|\F_m\right]E[W_i|\F_m]\right]\\
&=\; 0 ,\\
\end{split}
\end{equation}

where the second term vanishes by (\ref{eq:m_finite}).
In particular, for any indices $i,j$ with $i\neq j$ and any $W\in\L_p$,
we can let $W_i=W$, $W_j=-W$, and $W_\ell=0$ for $\ell\notin\{i,j\}$
to obtain $E[E[Y_i^\infty - Y_j^\infty|\F_m]W]=0$, from which it follows that
$E[Y_i^\infty|\F_m] = E[Y_j^\infty|\F_m]$ in $\L_q$; thus the sequence limit
$(Y_1^\infty,\ldots,Y_n^\infty,m)$ belongs to $K$.
\end{proof}

\section*{Acknowledgments}
This material is based upon work supported by the Massachusetts Clean Energy Center and the Maryland Energy Administrations as well as the U.S. Department of Energy’s Office of Energy Efficiency and Renewable Energy (EERE) under the Wind Energy Technologies Office (WETO) Award Number DE-EE0011269. The views expressed herein do not necessarily represent the views of the U.S. Department of Energy or the United States Government.

\bibliographystyle{elsarticle-num}
\bibliography{References}

\end{document}